\documentclass[reqno]{amsart}
\usepackage{amsmath, pb-diagram}
\usepackage{amssymb}
\usepackage{amscd}

\begin{document}
\title [ON FECKLY CLEAN RINGS]
 {ON FECKLY CLEAN RINGS}

\author{H. Chen}\thanks{\hspace{-7.0mm} $^1$Department of Mathematics, Hangzhou Normal University, Hangzhou 310034,
China, huanyinchen@aliyun.com}

\author{H. Kose}\thanks{\hspace{-6.0mm} $^2$Department of Mathematics, Ahi Evran University, Kirsehir,
Turkey, handankose@gmail.com}

\author{Y. Kurtulmaz}\thanks{\hspace{-6.0mm} $^3$Department of Mathematics, Bilkent University, Ankara, Turkey, yosum@fen.bilkent.edu.tr
}

\date{}
\newtheorem {thm}{Theorem}[section]
\newtheorem{lem}[thm]{Lemma}
\newtheorem{prop}[thm]{Proposition}
\newtheorem{cor}[thm]{Corollary}
\newtheorem{df}[thm]{Definition}
\newtheorem{nota}{Notation}
\newtheorem{note}[thm]{Remark}
\newtheorem{ex}[thm]{Example}
\newtheorem{exs}[thm]{Examples}
\newtheorem{rmk}[thm]{Remark}
\newtheorem{quo}[thm]{Question}

\def\bc{\begin{center}}
\def\ec{\end{center}}
\def\no{\noindent}
\def\hang{\hangindent\parindent}
\def\textindent#1{\indent\llap{[#1]\enspace}\ignorespaces}
\def\re{\par\hang\textindent}

\begin{abstract}
A ring $R$ is feckly clean provided that for any $a\in R$ there
exists an element $e\in R$ and a full element $u\in R$ such that
$a=e+u, eR(1-e)\subseteq J(R)$. We prove that a ring $R$ is feckly
clean if and only if for any $a\in R$, there exists an element
$e\in R$ such that $V(a)\subseteq V(e), V(1-a)\subseteq V(1-e)$
and $eR(1-e)\subseteq J(R)$, if and only if for any distinct
maximal ideals $M$ and $N$, there exists an element $e\in R$ such
that $e\in M, 1-e\in N$ and $eR(1-e)\subseteq J(R)$, if and only
if $J$-$spec(R)$ is strongly zero dimensional, if and only if
$Max(R)$ is strongly zero dimensional and every prime ideal
containing $J(R)$ is contained in a unique maximal ideal. More
explicit characterizations are also discussed for commutative
feckly clean rings.
\end{abstract}

\maketitle
\section{Introduction}

\vskip4mm An element $a$ of a ring $R$ is clean provided that $a$
is the sum of an idempotent and a unit in $R$. A ring $R$ is clean
provided that every element in $R$ is clean. Clean rings are
defined by Nicholson, and classes of rings that share properties
with clean rings are of great interest to many researchers (cf.
[1-3], [6], [8-10] and [14]). The motivation of this article is to
consider such a kind of rings and characterize them in terms of
topological properties. An element $u\in R$ is said to be a full
element if $RuR=R$. We say that an element $a\in R$ is feckly
clean provided that there exists an element $e\in R$ and a full
element $u\in R$ such that $a=e+u, eR(1-e)\subseteq J(R)$. A ring
$R$ is feckly clean provided that every element in $R$ is feckly
clean. It is easily shown that a ring $R$ is feckly clean if and
only if $R/J(R)$ is feckly clean. Clearly, every abelian clean
element is feckly clean. Thus, abelian clean rings are feckly
clean. But the converse is not true.

\vskip4mm \no {\bf Example 1.1.}\ \ Let $R=\{ \frac{m}{n}~|~m,n\in
{\Bbb Z}, (m,n)=1 ~\mbox{and}~ 2,3 \nmid n\}$. Then $Max(R)=\{
2R,3R\}$, and so $R/J(R)\cong {\Bbb Z}_2\oplus {\Bbb Z}_3$. Hence,
$R/J(R)$ is clean, and so $R$ is a commutative feckly clean ring.
As idempotents do not lift modulo $J(R)$, we see that $R$ is not
clean. In fact, $4\in R$ can not be written as the sum of an
idempotent and a unit.\\

Also right quasi-duo exchange rings are feckly clean. Some other
examples will be also provided. We will show that feckly clean
rings are mostly characterized by using the topological spaces of
all maximal ideals of $R$, and of all prime ideals containing the
Jacobson radical of $R$.

In Section 2, we show that a ring $R$ is feckly clean if and only
if for any $a\in R$, there exists an element $e\in R$ such that
$V(a)\subseteq V(e), V(1-a)\subseteq V(1-e)$ and $eR(1-e)\subseteq
J(R)$, if and only if for any distinct maximal ideals $M$ and $N$,
there exists an element $e\in R$ such that $e\in M, 1-e\in N$ and
$eR(1-e)\subseteq J(R)$. In Section 3, we consider feckly clean in
terms of the topological space of all prime ideals containing the
Jacobson radical of $R$. It is shown that a ring $R$ is feckly
clean if and only if $J$-$spec(R)$ is strongly zero dimensional,
if and only if $Max(R)$ is strongly zero dimensional and every
prime ideal containing $J(R)$ is contained in a unique maximal
ideal. In Section 4, we prove that every right (left) quasi-duo
ring is feckly clean, and provide more examples of such rings. In
the last section, more explicit characterizations are also
discussed for commutative feckly clean rings.

Throughout, all rings are associative with an identity. $Max(R)$
denote the denote the topological space of all maximal ideals of
$R$. We always use $J\mbox{-}spec(R)$ to stand for the topological
space of all prime ideals of $R$ that containing the Jacobson
radical. $J(R)$ will denote the Jacobson radical of $R$ and $U(R)$
will be the set of invertible element in $R$.

\section{Separative Properties}

\vskip4mm Let $R$ be a ring. Define the maximal spectrum of $R$ as
the set of all maximal ideals of $R$ which we denote as $Max(R)$.
Let $I$ be an ideal of $R$, and let $E(I)=\{P\in Max(R)\mid
I\not\subseteq P\}$. Then $E(R)=Max (R),E(0)=\emptyset, E(I)\cap
E(J)=E\big(IJ\big)$ and
$\bigcup\limits_{i}E(I_i)=E\big(\sum\limits_{i}I_i\big)$. So $Max
(R)$ is a topological space\index{space!topological}, where
$\{E(I)\mid I \unlhd R\}$ is the collection of its open sets. Let
$V(I)=Max(R)-E(I)$. Then $V(I)=\{P\in Max(R)\mid I\subseteq P\}$.
In addition, $V(0)=Max(R),V(R)=\emptyset, V(I)\bigcup
V(J)=V\big(IJ\big)$ and
$\bigcap\limits_{i}V(I_i)=V\big(\sum\limits_{i}I_i\big)$. Let
$a\in R$, and let $V(a):=V\big(RaR\big)$. It is easy to verify
that $V(I)=\bigcap\limits_{a\in I}V(a)$. Also we use $E(a)$ to
stand for $E\big(RaR\big)$. Let $S$ and $T$ be two sets. We always
use $S\sqcup T$ to denote the set $S\cup T$ with $S\cap
T=\emptyset$. In view of [4, Lemma 17.1.3], $Max(R)$ is a compact
space. The aim of this section is to characterize the feckly clean
rings by means of the separation property of the topological space
$Max(R)$.

\vskip4mm \no {\bf Lemma 2.1.}\ \ {\it Let $R$ be a ring, and let
$P$ be a maximal ideal of $R$. Then $J(R)\subseteq P$.} \vskip2mm
\no {\it Proof.}\ \ Let $x\in J(R)$. If $x\not\in P$, then
$RxR+P=R$. Write $s_1xt_1+\cdots +s_nxt_n+y=1$ where $y\in P,
s_i,t_i\in R$ for $1\leq i\leq n$. Then $y=1-(s_1xt_1+\cdots
+s_nxt_n)\in U(R)$, a contradiction. Therefore $x\in P$.
Accordingly, $J(R)\subseteq P$.\hfill$\Box$

\vskip4mm \no {\bf Lemma 2.2.}\ \ {\it Let $R$ be a ring. Then the
following are equivalent:} \vspace{-.5mm}
\begin{enumerate}
\item [(1)]{\it $a\in R$ is feckly clean.} \vspace{-.5mm}
\item [(2)]{\it There exists an
element $e\in R$ such that $V(a-1)\subseteq V(e)\subseteq
Max(R)-V(a), eR(1-e)\subseteq J(R)$.}
\end{enumerate}
\vspace{-.5mm}  {\it Proof.}\ \ $(1)\Rightarrow (2)$ Let $a\in R$
be feckly clean. Then there exist an $e\in R$ and a full element
$u\in U(R)$ such that $a=e+u$ and $eR(1-e)\subseteq J(R)$. For any
$P\in V(a-1)$, we see that $P\not\in V(1-e)$; otherwise,
$u=(a-1)+(1-e)\in P$. Clearly, $V(e)\bigcap V(1-e)=\emptyset$. For
any $P\in Max(R)$, as $eR(1-e)\subseteq J(R)$, we see that
$eR(1-e)\subseteq P$ by Lemma 2.1. This implies that $e\in P$ or
$1-e\in P$. Thus, $P\in V(e)$ or $P\in V(1-e)$. Therefore
$Max(R)=V(e)\bigsqcup V(1-e)$; hence, $P\in V(e)$. We infer that
$V(a-1)\subseteq V(e)$. If $P\in Max(R)$ and $P\not\in
Max(R)-V(a)$, then $P\in V(a)$; hence, $P\not\in V(e)$; otherwise,
$u=a-e\in P$. This implies that $V(e)\subseteq Max(R)-V(a).$
Therefore $V(a-1)\subseteq V(e)\subseteq Max(R)-V(a)$.

$(2)\Rightarrow (1)$ Assume that there exists an element $e\in R$
such that $V(a-1)\subseteq V(e)\subseteq Max(R)-V(a),
eR(1-e)\subseteq J(R)$ Let $u=a-e$. If $RuR\neq R$, by using
Zorn's Lemma, there exists a maximal ideal $P$ such that
$RuR\subseteq P\subsetneqq R$. In light of Lemma 2.1,
$eR(1-e)\subseteq J(R)\subseteq P$, we see that $e\in P$ or
$1-e\in P$. Hence, $P\in V(e)$ or $P\in V(1-e)$. If $P\in V(e)$,
then $a=e+u\in P$, whence, $P\in V(a)$, a contradiction. If $P\in
V(1-e)$, then $a-1=(e-1)+u\in P$, whence, $P\in V(a-1)$. Thus
$P\in V(e)$, a contradiction. Therefore $a-e\in R$ is a full
element, and so $a$ is feckly clean.\hfill$\Box$

\vskip4mm \no {\bf Theorem 2.3.}\ \ {\it Let $R$ be a ring. Then
the following are equivalent:} \vspace{-.5mm}
\begin{enumerate}
\item [(1)]{\it $R$ is feckly clean.}
\vspace{-.5mm}
\item [(2)]{\it For any
disjoint closed sets $A$ and $B$ of $Max(R)$, there exists an
element $e\in R$ such that $A\subseteq V(e), B\subseteq V(1-e)$
and $eR(1-e)\subseteq J(R)$.} \vspace{-.5mm}
\item [(3)]{\it For any $a\in R$, there exists an
element $e\in R$ such that $V(a)\subseteq V(e), V(1-a)\subseteq
V(1-e)$ and $eR(1-e)\subseteq J(R)$.}
\end{enumerate} \vspace{-.5mm} {\it Proof.}\ \ $(1)\Rightarrow (2)$ Let $A$ and $B$ be disjoint closed sets of
$Max(R)$. Then $A\bigcap B=\emptyset$. Clearly, there exist two
ideals $I$ and $J$ such that $A=V(I)$ and $B=V(J)$; hence,
$V(I)\bigcap V(J)=\emptyset$. If $I+J\neq R$, then there exists a
maximal ideal $P$ of $R$ such that $I+J\subseteq P\subsetneqq R$.
Hence, $P\in V(I+J)=V(I)\bigcap V(J)$, a contradiction. This
implies that $I+J=R$, and so $a+b=1$ for some $a\in I$ and $b\in
J$. As $a\in R$ is feckly clean, it follows from Lemma 2.2 that
there exists an element $e\in R$ such that
$$V(a-1)\subseteq V(1-e)\subseteq Max(R)-V(a), (1-e)Re\in J(R)$$   It is easy to check
that
$$\begin{array}{l}
B= V(J)\subseteq V(b)\\
=V(a-1)\subseteq V(1-e)\subseteq Max(R)-V(a)\\
\subseteq Max(R)-V(I)= Max(R)-A.\end{array}$$ Thus $B\subseteq
V(1-e)$. As $Max(R)=V(e)\bigsqcup V(1-e)$ and $V(1-e)\subseteq
Max(R)-A$, we see that $A\subseteq V(e)$, as required.

$(2)\Rightarrow (3)$ For any $a\in R$, $V(a)\bigcap
V(1-a)=\emptyset$. That is, $V(a)$ and $V(a-1)$ are disjoint
closed sets of $Max(R)$. By hypothesis, there exists an element
$e\in R$ such that $V(a)\subseteq V(e), V(1-a)\subseteq V(1-e)$
and $eR(1-e)\subseteq J(R)$.

$(3)\Rightarrow (1)$ For any $a\in R$, we can find an element
$e\in R$ such that $V(a)\subseteq V(1-e)$ and $V(a-1)\subseteq
V(e)$ with $(1-e)Re\subseteq J(R)$. If $P\in Max(R)$ and $P\not\in
Max(R)-V(a)$, then $P\in V(a)$, and so $P\in V(1-e)$. This implies
that $P\not\in V(e)$. Hence, $V(e)\subseteq Max(R)-V(a)$. Clearly,
$\big(ReR(1-e)R\big)^2\subseteq R(1-e)ReR\subseteq J(R)$. As
$J(R)$ is semiprime, we see that $eR(1-e)\subseteq J(R)$. In view
of Lemma 2.2, $a\in R$ is feckly clean.\hfill$\Box$

\vskip4mm \no  {\bf Corollary 2.4.}\ \ {\it Let $R$ be a ring.
Then the following are equivalent:} \vspace{-.5mm}
\begin{enumerate}
\item [(1)]{\it $R$ is feckly clean.}
\vspace{-.5mm}
\item [(2)]{\it For any disjoint
compact sets $A$ and $B$ of $Max(R)$, there exists an element
$e\in R$ such that $A\subseteq V(e), B\subseteq V(1-e)$ and
$eR(1-e)\subseteq J(R)$.} \vspace{-.5mm}
\end{enumerate} \vspace{-.5mm} {\it Proof.}\ \ $(1)\Rightarrow
(2)$ Let $P,Q\in Max(R), P\neq Q$. Then there exist $a\in P, b\in
Q$ such that $a+b=1$. As $R$ is feckly clean, we can find an
element $e\in R$ and a full element $u\in R$ such that $a=e+u$ and
$eR(1-e)\subseteq J(R)$. This implies that $e\not\in P$. Clearly,
$b=1-a=(1-e)-u$, and so $1-e\not\in Q$. Consequently, $1-e\in P$
and $e\in Q$. Thus, $P\in V(1-e)$ and $Q\in V(e)$. One easily
checks that $Max(R)=V(e)\bigsqcup V(1-e)$. Hence, $V(e)$ and
$V(1-e)$ are clopen sets. That is, $Max(R)$ is an Hausdorff space.
In light of [4, Lemma 17.1.3], $Max(R)$ is compact. So every
compact subset of $Max(R)$ is closed. For any disjoint compact
subsets $A$ and $B$ of $Max(R)$, it follows from Lemma 2.2 that
there exists an element $e\in R$ such that $A\subseteq V(e),
B\subseteq V(1-e)$ and $eR(1-e)\subseteq J(R)$.

$(2)\Rightarrow (1)$ As $Max(R)$ is compact, every closed set is
compact. Thus we complete the proof by Theorem 2.3.\hfill$\Box$

\vskip4mm \no {\bf Lemma 2.5.}\ \ {\it If $e,f\in R$ such that
$eR(1-e), fR(1-f)\subseteq J(R)$, then there exists $g\in R$ such
that $E(e)\bigcup E(f)=E(g)$ and $gR(1-g)\subseteq J(R)$.}
\vskip2mm\no {\it Proof.}\ \ Set $g=e+f-ef$. If $P\not\in
E(e)\bigcup E(f)$, then $e, f\in P$, and so $e+f-ef\in P$. This
implies that $P\not\in E(e+f-ef)$; hence, $E(g)\subseteq
E(e)\bigcup E(f)$. Suppose that $P\in E(e)\bigcup E(f)$, then
$e\not\in P$ or $f\not\in P$. As $eR(1-e), fR(1-f)\subseteq
J(R)\subseteq P$, we see that $1-e\in P$ or $1-f\in P$. This
implies that $(1-e)(1-f)\in P$, and then $1-(1-e)(1-f)\not\in P$.
That is, $e+f-ef\not\in P$. Therefore $P\in E(g)$, and so
$E(e)\bigcup E(f)\subseteq E(g)$. As a result, $E(e)\bigcup
E(f)=E(g)$. On the other hand, $gR(1-g)\subseteq
eR(1-e)(1-f)+(1-e)fR(1-e)(1-f)\subseteq J(R)$, as asserted.
\hfill$\Box$

\vskip4mm \no {\bf Theorem 2.6.}\ \ {\it Let $R$ be a ring. Then
the following are equivalent:} \vspace{-.5mm}
\begin{enumerate}
\item [(1)]{\it $R$ is feckly clean.}
\vspace{-.5mm}
\item [(2)]{\it For any distinct maximal
ideals $M$ and $N$, there exists an element $e\in R$ such that
$e\in M, 1-e\in N$ and $eR(1-e)\subseteq J(R)$.}
\end{enumerate} \vspace{-.5mm} {\it Proof.}\ \ $(1)\Rightarrow
(2)$ Suppose that $M$ and $N$ are distinct maximal ideals. Then
$M+N=R$. Write $1=a+b$ with $a\in M, b\in N$. As $R$ is feckly
clean, there exist an element $e\in R$ and a full element $u\in R$
such that $b=e+u$ and $eR(1-e)\subseteq J(R)$. Obviously,
$Max(R)=V(e)\bigsqcup V(1-e)$. If $M\in V(1-e)$, then $1-e\in M$.
As $1-b=a\in M$, we see that $u\in M$, a contradiction. Therefore
$M\in V(e)$, i.e., $e\in M$. Similarly, we show that $1-e\in N$,
as desired.

$(2)\Rightarrow (1)$ Let $A$ and $B$ be disjoint closed sets of
$Max(R)$. Take $M\in B$. For any $N\in A$, $M\neq N$. By
hypothesis, there exists an element $g\in R$ such that $g\in M,
1-g\in N$ and $gR(1-g)\subseteq J(R)$. Hence, $g\not\in N$, and so
$N\in E(g)$. As $Max(R)$ is compact, and then so is $A$. Thus, we
can find some $g_1,\cdots ,g_n\in R$ such that $A\subseteq
\bigcup\limits_{i=1}^{n}E(g_i)$ and $g_iR(1-g_i)\subseteq J(R)$
for each $i$. According to Lemma 2.5, we can find an element $f\in
R$ such that $\bigcup\limits_{i=1}^{n}E(g_i)=E(f)$ and
$fR(1-f)\subseteq J(R)$. Set $e=1-f$. Then $(1-e)Re\subseteq
J(R)$. This implies that $eR(1-e)\subseteq J(R)$. Further, we see
that $A\subseteq V(e)$.

Clearly, $M\in V(g_1)$ and $M\in V(g_2)$. Hence, $M\in
V(g_1+g_2-g_1g_2)$. By iteration of this process, $M\in
V(f)=V(1-e)$. As $B$ is compact, we can find some $e_1,\cdots
,e_m\in R$ such that $A\subseteq V(e_i), B\subseteq
\bigcup\limits_{i=1}^{m}E(e_i)$ and $e_iR(1-e_i)\subseteq J(R)$
for each $i$. In light of Lemma 2.5, we can find an element
$\alpha\in R$ such that $B\subseteq E(\alpha)$. Clearly,
$V(e_1)\bigcap V(e_2)\subseteq V(e_1+e_2-e_1e_2)$. Furthermore,
$A\subseteq \bigcap\limits_{i=1}^{m}V(e_i)=V(\alpha)$. In
addition, $\alpha R(1-\alpha)\subseteq J(R)$. According to Theorem
2.3, $R$ is feckly clean.\hfill$\Box$

\vskip4mm \no {\bf Corollary 2.7.}\ \ {\it If $a+b =1$ in $R$
implies that there exists an element $e\in R$ such that $1+ar\in
eR, 1+bs\in (1-e)R$ where $eR(1-e)\subseteq J(R)$, then $R$ is
feckly clean.} \vskip2mm \no {\it Proof.}\ \ Assume that $M$ and
$N$ are distinct maximal ideals of $R$. Then $M+N=R$. Write
$a+b=1$ with $a\in M, b\in N$. By hypothesis, there exist $r,s\in
R$ and an element $e\in R$ such that $1+ar\in eR, 1+bs\in (1-e)R$,
where $eR(1-e)\subseteq J(R)$. As $a\in M$, we see that $e\not\in
M$. As $eR(1-e)\subseteq J(R)\subseteq M$, we have $1-e\in M$.
Similarly, $1-e\not\in N$, and so $e\in N$. In light of Theorem
2.6, $R$ is feckly clean.\hfill$\Box$

\vskip4mm A ring $R$ is local if it has only one maximal right
ideal. As is well known, a ring $R$ is local if and only if
$a+b=1$ in $R$ implies that either $a$ or $b$ is invertible. We
claim that every local ring $R$ is feckly clean. Given $a+b=1$ in
$R$. Then either $a\in U(R)$ or $1-a\in U(R)$. If $a\in U(R)$,
choosing $e=0$, then $1+a(-a^{-1})\in eR, 1+b\cdot 0\in (1-e)R$.
If $b\in U(R)$, choosing $e=1$, then $1+a\cdot 0\in eR, 1+b\cdot
(-b^{-1})\in (1-e)R$. Clearly, $eR(1-e)\subseteq J(R)$ in each
case, and we are done.

\section{Strongly Zero-dimensional Spaces}

\vskip4mm A topological space $X$ is said to be strongly
zero-dimensional\index{space!strongly zero-dimensional} provided
that any two disjoint closed sets are separated by clopen
sets\index{set!clopen}, that is, if $A$ and $B$ are disjoint
closed sets\index{set!closed}, then there exist disjoint clopen
sets $C_1,C_2$ such that $A\subseteq C_1$ and $B\subseteq C_2$.
The aim of this section is to investigate feckly clean rings by
means of some strongly zero-dimensional spaces.

Let $J$-$spec(R)=\{ P\in Spec(R)~|~J(R)\subseteq P\}$. As is well
known, the Jacobson radical $J(R)$ is semiprime, and so $J(R)$ is
the intersection of some prime ideals. Hence,
$J(R)=\bigcap\limits_{P\in J-\mbox{Spec}(R)}P$. Let $I$ be an
ideal of $R$, and let $F(I)=\{P\in J\mbox{-}spec(R)\mid
I\not\subseteq P\}$. Then $F(R)=J\mbox{-}spec(R), F(0)=\emptyset,
F(I)\cap F(J)=F\big(IJ\big)$ and
$\bigcup\limits_{i}F(I_i)=F\big(\sum\limits_{i}I_i\big)$. So
$J$-$spec(R)$ is a topological space, where $\{F(I)\mid I \unlhd
R\}$ is the collection of its open sets. Let
$W(I)=J\mbox{-}spec(R)-F(I)$. Then $W(I)=\{P\in
J\mbox{-}spec(R)\mid I\subseteq P\}$ is the collection of its
closed sets. We use $W(a)$ to denote $W\big(RaR\big)$ for any
$a\in R$.

\vskip4mm \no {\bf Lemma 3.1.}\ \ {\it Let $R$ be a ring. Then
$a\in R$ is feckly clean if and only if there exists an element
$e\in R$ such that $W(a-1)\subseteq W(e)\subseteq
J\mbox{-}spec(R)-W(a), eR(1-e)\subseteq J(R)$.}\vskip2mm \no {\it
Proof.}\ \ Let $a\in R$ be feckly clean. Then there exist an
element $e\in R$ and a full element $u\in R$ such that $a=e+u$ and
$eR(1-e)\subseteq J(R)$. For any $P\in J\mbox{-}spec(R)$, as
$eR(1-e)\subseteq J(R)\subseteq P$, we see that $e\in P$ or
$1-e\in P$. Thus, $P\in W(e)$ or $P\in W(1-e)$. Therefore
$J\mbox{-}spec(R)=W(e)\bigsqcup W(1-e)$. As in the proof of Lemma
2.2, $W(a-1)\subseteq W(e)\subseteq J\mbox{-}spec(R)-W(a)$, as
required.

The proof of the converse is the same as in Lemma 2.2.\hfill$\Box$

\vskip4mm \no {\bf Lemma 3.2.}\ \ {\it Let $R$ be a ring. If $A$
is a clopen subset of $J\mbox{-}spec(R)$, then there exists an
element $e\in R$ such that $A=W(e), eR(1-e)\subseteq J(R)$.}
\vskip2mm \no {\it Proof.}\ \ Let $A$ be a clopen subset of
$J\mbox{-}spec(R)$. Then we can find ideals $I$ and $J$ of $R$
such that $A=W(I)$ and $J\mbox{-}spec(R)=W(I)\sqcup W(J)$. Thus,
$W(IJ)=J\mbox{-}spec(R)$. Hence $IJ\subseteq \bigcap\limits_{P\in
Spec(R)}\{ P~\mid~J(R)\subseteq P\}=J(R)$. On the other hand,
$W(I)\bigcap W(J)=\emptyset$, and so $W(I+J)=\emptyset$. If
$I+J\neq R$, then there exists a maximal ideal $P$ of $R$ such
that $I+J\subseteq P\subsetneqq R$. In light of Lemma 2.1,
$J(R)\subseteq P$, and so $P\in J\mbox{-}spec(R)$; hence, $P\in
W(I+J)$, a contradiction. Thus, $I+J=R$. Write $1=e+f, e\in I,f\in
J$. Then
$$J\mbox{-}spec(R)=W(e)\sqcup W(f), W(I)\subseteq W(e), W(J)\subseteq
W(f).$$ It follows from $J\mbox{-}spec(R)=W(I)\sqcup W(J)$ that
$W(I)=W(e)$, as desired.\hfill$\Box$

\vskip4mm \no {\bf Theorem 3.3.}\ \ {\it A ring $R$ is feckly
clean if and only if $J$-$spec(R)$ is strongly zero
dimensional.}\vskip2mm \no {\it Proof.}\ \ Suppose that $R$ is
feckly clean. Let $A$ and $B$ be disjoint closed sets of
$J$-$spec(R)$. Then $A\bigcap B=\emptyset$. Clearly, there exist
two ideals $I$ and $J$ such that $A=W(I)$ and $B=W(J)$. It is easy
to verify that $I+J=R$. Write $a+b=1$ with $a\in I$ and $b\in J$.
Since $a\in R$ is feckly clean, by virtue of Lemma 3.1, there
exists an element $e\in R$ such that
$$W(a-1)\subseteq W(1-e)\subseteq J\mbox{-}spec(R)-W(a), (1-e)Re\subseteq J(R).$$  Clearly,
$$\begin{array}{l}
B= W(J)\subseteq W(b)\\
=W(a-1)\subseteq W(1-e)\subseteq J\mbox{-}spec(R)-W(a)\\
\subseteq J\mbox{-}spec(R)-W(I)= J\mbox{-}spec(R)-A.\end{array}$$
Obviously, $B\subseteq W(1-e)$. Since $W(1-e)\subseteq
J\mbox{-}spec(R)-A$, we get $A\subseteq W(e)$. In addition,
$eR(1-e)\subseteq J(R).$ Therefore $J$-$spec(R)$ is strongly zero
dimensional.

Conversely, assume that $J$-$spec(R)$ is strongly zero
dimensional. For any $a\in R$, we see that $W(a)$ and $W(a-1)$ are
disjoint closed sets of $J\mbox{-}spec(R)$. Thus, we can find two
distinct clopen sets $U$ and $V$ such that $W(a)\subseteq U$ and
$W(a-1)\subseteq V$. In light of Lemma 3.2, we have elements $e,
f\in R$ such that $U=W(e)$ and $V=W(f)$, where $eR(1-e),
fR(1-f)\subseteq J(R)$. If $P\in J\mbox{-}spec(R)$ and $P\not\in
J\mbox{-}spec(R)-W(a)$, then $P\in W(a)$, and so $P\in W(e)$.
Since $W(e)\bigcap W(f)=W(ReR+RfR)=\emptyset$, it is easy to check
that $ReR+RfR=R$. Write
$\sum\limits_{i=1}^{m}s_iet_i+\sum\limits_{j=1}^{n}x_jfy_j=1$. If
$P\in W(f)$, then $RfR\subseteq P$, and so
$\sum\limits_{i=1}^{m}s_iet_i\not\in P$. This implies that
$e\not\in P$. As $eR(1-e)\subseteq J(R)\subseteq P$, we see that
$1-e\in P$. This shows that $P\in W(1-e)$, a contradiction. Thus,
$P\not\in W(f)$. As a result,
$$W(a-1)\subseteq W(f)\subseteq J\mbox{-}spec(R)-W(a).$$ In view of Lemma 3.1, $a\in R$ is feckly clean. As a result, $R$ is feckly clean.\hfill$\Box$

\vskip4mm Recall that a ring $R$ is a Hilbert ring if every prime
ideal of the ring is an intersection of maximal ideals (cf. [6]).
For instance, any artinian ring, ring of integers, any polynomial
ring in finitely many variables over a field. A general form
states that if $R$ is a Hilbert ring, then so is any finitely
generated $R$-algebra $S$.

\vskip4mm \no {\bf Corollary 3.4.}\ \ {\it Let $R$ be a Hilbert
ring. Then $R$ is feckly clean if and only if $Max(R)$ is strongly
zero-dimensional.}\vskip2mm \no {\it Proof.}\ \ By hypothesis,
every prime ideal contains the Jacobson radical as in the proof of
Lemma 2.1. Therefore the result follows from Theorem
3.3.\hfill$\Box$

\vskip4mm \no {\bf Theorem 3.5.}\ \ {\it A ring $R$ is feckly
clean if and only if } \vspace{-.5mm}
\begin{enumerate}
\item [(1)]{\it $Max(R)$ is strongly zero dimensional.}
\vspace{-.5mm}
\item [(2)]{\it Every prime ideal of $R$ containing $J(R)$ is contained in a unique maximal ideal.}
\end{enumerate} \vspace{-.5mm} {\it Proof.}\ \ Suppose that $R$ is feckly clean. In view of Theorem 2.3,
$Max(R)$ is strongly zero dimensional. If there exists an ideal
$P\in J\mbox{-}spec(R)$ such that $P\subseteq M\bigcap N$ where
$M$ and $N$ are distinct maximal ideals, by virtue of Theorem 2.6,
there exists an element $e\in R$ such that $e\in M, 1-e\in N$ and
$eR(1-e)\subseteq J(R)$. As $eR(1-e)\subseteq J(R)\subseteq P$,
$e\in P$ or $1-e\in P$. If $e\in P$, then $e\in N$, a
contradiction. If $1-e\in P$, then $1-e\in M$, a contradiction.
Therefore every prime ideal of $R$ containing $J(R)$ is contained
in exactly one maximal ideal.

Conversely, assume $(1)$ and $(2)$ hold. Then there exists a map
$\varphi: J\mbox{-}spec(R)\to Max(R)$, $\varphi(P)=M$, where $M$
is the unique maximal ideal such that $P\subseteq M$. Let $I$ be
an ideal of $R$, we let $V_S(I)=\{ P\in
J\mbox{-}spec(R)~|~I\subseteq P\}$ and $V_M(I):=\{ P\in
Max(R)~|~I\subseteq P\}$. Then $\varphi\big(V_S(I)\big)=V_M(I)$.
For any disjoint closed sets $A, B\subseteq J\mbox{-}spec(R)$,
there exist two ideals $I$ and $J$ of $R$ such that $A=V_S(I)$ and
$B=V_S(J)$. Thus, $\varphi(A)$ and $\varphi(B)$ are both closed.
As $V_S(I)\bigcap V_S(J)=\emptyset$, $V_S(I+J)=\emptyset$; hence,
$I+J=R$. This implies that $V_M(I)\bigcap
V_M(J)=V_M(I+J)=V_M(R)=\emptyset$. That is, $\varphi(A)$ and
$\varphi(B)$ are disjoint closed sets of $Max(R)$. By hypothesis,
$Max(R)$ is strongly zero-dimensional, and so we can find disjoint
clopen sets $U, V\subseteq Max(R)$ such that $V_M(I)\subseteq U,
V_M(J)\subseteq V$. Clearly, $A\subseteq \varphi^\leftarrow
\varphi (A)\subseteq \varphi^\leftarrow (U)$ and $B\subseteq
\varphi^\leftarrow \varphi (B)\subseteq \varphi^\leftarrow (V)$.
As in the proof of [7, Theorem 1.2], $\varphi$ is continuous;
hence, $\varphi^\leftarrow (U)$ and $\varphi^\leftarrow (V)$ are
clopen. For any $P\in \varphi^\leftarrow (U)\bigcap
\varphi^\leftarrow (V)$, there exists a unique $M\in Max(R)$ such
that $P\subseteq M$. Hence, $M\in U\bigcap V$, a contradiction.
This implies that $\varphi^\leftarrow (U)\bigcap
\varphi^\leftarrow (V)=\emptyset$. Consequently,
$J\mbox{-}spec(R)$ is strongly zero-dimensional. In light of
Theorem 3.3, we complete the proof.\hfill$\Box$

\vskip4mm \no {\bf Corollary 3.6.}\ \ {\it A ring $R$ is feckly
clean if and only if } \vspace{-.5mm}
\begin{enumerate}
\item [(1)]{\it $Max(R)$ is strongly zero dimensional.}
\vspace{-.5mm}
\item [(2)]{\it For any distinct maximal
ideals $M$ and $N$, there exist elements $a, b\in R$ such that
$a\not\in M, b\not\in N$ and $aRb\subseteq J(R)$.}
\end{enumerate} \vspace{-.5mm} {\it Proof.}\ \ Suppose that $R$ is feckly
clean. Then $Max(R)$ is strongly zero dimensional by Theorem 3.5.
For any distinct maximal ideals $M$ and $N$, it follows from
Theorem 2.6 that there exists an element $e\in R$ such that $e\in
M, 1-e\in N$ and $eR(1-e)\subseteq J(R)$. Set $a=1-e$ and $b=e$.
Then $a\not\in M, b\not\in N$ and $aRb=(1-e)Re\subseteq J(R)$.

$(2)\Rightarrow (1)$ Suppose that $P$ containing $J(R)$ is a prime
ideal such that $P\subseteq M\bigcap N$ where $M$ and $N$ are
distinct maximal ideals of $R$. Then $aRb\subseteq J(R)\subseteq
P$; hence, either $a\in P$ or $b\in P$. This implies that $a\in M$
or $b\in N$, a contradiction. Therefore every prime ideal of $R$
containing $J(R)$ is contained in a unique maximal ideal.
According to Theorem 3.5, we complete the proof.\hfill$\Box$

\vskip4mm \no {\bf Corollary 3.7.}\ \ {\it If $a+b=1$ in $R$
implies that $(1+ar)R(1+bs)\subseteq J(R)$ for some $r,s\in R,$
then the following are equivalent:}\vspace{-.5mm}
\begin{enumerate}
\item [(1)]{\it $R$ is feckly clean.} \vspace{-.5mm}
\item [(2)]{\it $Max(R)$ is strongly zero dimensional.}
\end{enumerate} \vspace{-.5mm} {\it Proof.}\ \ $(1)\Rightarrow (2)$ In view of Theorem 3.5,
$Max(R)$ is strongly zero dimensional.

$(2)\Rightarrow (1)$ Let $P\in J\mbox{-}spec(R)$. Then $P$ is
contained in some maximal ideal of $R$. Suppose that $P$ is
contained in two distinct maximal ideals $M$ and $N$. Then
$M+N=R$. Write $a+b=1$ with $a\in M, b\in N$. By hypothesis, there
are some $r,s\in R$ such that $(1+ar)R(1+bs)\subseteq J(R)$. As
$J(R)\subseteq P$, we see that either $1+ar\in P$ or $1+bs\in P$.
If $1+ar\in P$, then $1=(1+ar)-ar\in M$. If $1+bs\in P$, then
$1=(1+bs)-bs\in N$. In any case, we get a contradiction. Therefore
every prime ideal of $R$ containing $J(R)$ is contained in exactly
one maximal ideal. According to Theorem 3.5, $R$ is feckly
clean.\hfill$\Box$

\section{More Examples}

\vskip4mm Recall that a ring $R$ is $\pi$-regular provided that
for any $a\in R$ there exists $n\in {\Bbb N}$ such that $a^n\in
a^nRa^n$.

\vskip4mm \no {\bf Theorem 4.1.}\ \ {\it If $R/J(R)$ is an abelian
$\pi$-regular ring, then $R$ is feckly clean.} \vskip2mm\no {\it
Proof.}\ \ Let $a\in R$. In view of [13, Theorem 30.2], $R/J(R)$
is clean. Thus, there exists an idempotent $\overline{e}\in
R/J(R)$ and a unit $\overline{u}\in R/J(R)$ such that
$\overline{a}=\overline{e}+\overline{u}$. As units lift modulo
$J(R)$, we may assume that $u\in U(R)$. Thus, we can find $r\in
J(R)$ such that $a=e+(u+r)$. Clearly, $u+r\in U(R)$, and so
$u+r\in R$ is a full element. As $R/J(R)$ is abelian, we see that
$\overline{e}\big(R/J(R)\big)(\overline{1}-\overline{e})=\overline{0}$.
Therefore $eR(1-e)\subseteq J(R)$. Therefore $R$ is feckly clean.
\hfill$\Box$

\vskip4mm Sz$\acute{a}$sz in [12] studied the ring $R$ with the
property that for any $x\in R$, there exists some $n(x)\geq 2$
such that $xRx=x^{n(x)}Rx^{n(x)}$. He called these rings
$gsr$-rings. For instance, every strongly regular ring is a
$gsr$-ring.

\vskip4mm \no {\bf Corollary 4.2.}\ \ {\it Every $gsr$-ring is
feckly clean.}\vskip2mm \no {\it Proof.}\ \ Let $R$ be a
$gsr$-ring. Then for any $x\in R$, $x^2\in x^2Rx^2$; hence, $R$ is
a $\pi$-regular ring. This implies that $R/J(R)$ is $\pi$-regular.
Given $\overline{x}^2=\overline{0}$ in $R/J(R)$, then $x^2\in
J(R)$. By hypothesis, $xRx=x^2Rx^2\subseteq J(R)$, i.e.,
$(RxR)^2\subseteq J(R)$. As $J(R)$ is semiprime, it follows that
$RxR\subseteq J(R)$, and so $x\in J(R)$. That is,
$\overline{x}=\overline{0}$. This implies that $R/J(R)$ is
reduced. For any idempotent $e\in R/J(R)$ and any $a\in R/J(R)$,
it follows from $\big(ea(\overline{1}-e)\big)^2=0$ that
$ea(\overline{1}-e)=0$, thus $ea=eae$. Likewise, $ae=eae$. This
implies that $ea=ae$. As a result, every idempotent in $R/J(R)$ is
central, i.e., $R/J(R)$ is abelian. Therefore we complete the
proof by Theorem 4.1.\hfill$\Box$

\vskip4mm \no {\bf Corollary 4.3.}\ \ {\it If for any $a\in R$,
there exists some $n\geq 2$ such that $a^n-a\in J(R)$, then $R$ is
feckly clean.} \vskip2mm\no {\it Proof.}\ \ Let $a\in R$. Then
there exists some $n\geq 2$ such that $a^n-a\in J(R)$, and so
$\overline{a}=\overline{a^n}$. Therefore $R/J(R)$ is regular. By
the proof of Corollary 4.2, $R/J(R)$ is abelian. Hence $R$ is
feckly clean by Theorem 4.1\hfill$\Box$

\vskip4mm \no {\bf Corollary 4.4.}\ \ {\it If $R/J(R)$ is a finite
commutative ring, then $R$ is feckly clean.} \vskip2mm\no {\it
Proof.}\ \ As $R/J(R)$ is a finite commutative ring, for any $a\in
R$, there exist distinct $m,n\in {\Bbb N}$ such that
$\overline{a^m}=\overline{a^n}$. Therefore $R/J(R)$ is
$\pi$-regular, and then the result follows by Theorem
4.1.\hfill$\Box$

\vskip4mm A ring $R$ a right (left) quasi-duo ring if every
maximal right (left) ideal of $R$ is an ideal. For instance, local
rings, duo rings and weakly right (left) duo rings are all right
(left) quasi-duo rings (cf. [15]). Further, every abelian exchange
ring is a right (left) quasi-duo ring, where a ring $R$ is an
exchange ring provided that for any $a\in R$ there exists an
idempotent $e\in R$ such that $e\in aR$ and $1-e\in (1-a)R$ (cf.
[4] and [15]).

\vskip4mm \no {\bf Theorem 4.5.}\ \ {\it Every right (left)
quasi-duo exchange ring is feckly clean.} \vskip2mm\no {\it
Proof.}\ \ Let $R$ be a right quasi-duo exchange ring. In view of
[15, Lemma 2.3], $R/J(R)$ is abelian. As $R$ is an exchange ring,
then so is $R/J(R)$. In light of [4, Theorem 17.2.2], $R/J(R)$ is
clean. Let $x\in R$. Then there exists an idempotent
$\overline{e}\in R/J(R)$ and a unit $\overline{u}\in R/J(R)$ such
that $\overline{a}=\overline{e}+\overline{u}$. We may assume that
$u\in U(R)$. Thus, we can find a $r\in J(R)$ such that
$a=e+(u+r)$, where $u+r\in U(R)$ is a full element. As $R/J(R)$ is
abelian, we see that $eR(1-e)\subseteq J(R)$, and so $R$ is feckly
clean. \hfill$\Box$

\vskip4mm \no {\bf Remark 4.6.}\ \ Let $R=\left(
\begin{array}{cc}
{\Bbb Z}_4&{\Bbb Z}_4\\
0&{\Bbb Z}_4 \end{array} \right)$. Then $R$ is a quasi-duo
exchange ring (cf. [15]). According to Theorem 4.5, $R$ is feckly
clean. Note that $R$ is not abelian.

\vskip4mm \no{\bf Remark 4.7.}\ \ Let $p$ and $q$ be two distinct
primes other than $2$. Then the ring ${\Bbb Z}_{(p)}\bigcap {\Bbb
Z}_{(q)}$ is feckly clean, but it is not clean. As in the proof of
[1, Proposition 16], we see that ${\Bbb Z}_{(p)}\bigcap {\Bbb
Z}_{(q)}$ is commutative with exactly two maximal ideals.
Additionally, $J(R)=pR\bigcap qR$. Therefore $R$ is a semilocal
ring. That is, $R/J(R)$ is finite direct sums of division rings.
Hence, $R/J(R)$ is clean. Consequently, $R$ is feckly clean. Since
$(p,q)=1$, $\frac{p(q+1)}{p+q}\in R$. Observing that $R$ is an
integral domain, the set of all idempotents in $R$ is
$\{\frac{0}{1}, \frac{1}{1}\}$. As $q\nmid q+1$, we see that
$\frac{p(q+1)}{p+q}\not\in U(R)$. As $p\nmid (p-1)q$, we see that
$\frac{p(q+1)}{p+q}-\frac{1}{1}\not\in U(R)$. This shows that
$\frac{p(q+1)}{p+q}\in R$ is not clean.

\section{The Commutative Case}

\vskip4mm If $R$ is commutative then $R$ is feckly clean if and
only if $R/J(R)$ is clean. In other words, in the commutative
case, the difference between "feckly clean" and "clean" is that
both are clean modulo the radical but in the former idempotents do
not necessarily lift modulo the radical. The aim of this section
is to investigate the necessary and sufficient conditions under
which a commutative ring is feckly clean.

\vskip4mm \no {\bf Theorem 5.1.}\ \ {\it Let $R$ be a commutative
ring. Then the following are equivalent:} \vspace{-.5mm}
\begin{enumerate}
\item [(1)]{\it $R$ is feckly clean }
\item [(2)]{\it For any $a\in R$, there exists $e\in R$ such that $a-e\in U(R), e^2-e\in J(R)$.}
\item [(3)]{\it For any $a\in R$, there exists $e\in R$ such that
$a-e\in (a-a^2)R, e-e^2\in J(R)$.}
\item [(4)]{\it For any $a\in R$, there exists $e\in R$ such that
$e\in aR, 1-e\in (1-a)R, e-e^2\in J(R)$.}
\end{enumerate}
\vspace{-.5mm}  {\it Proof.}\ \ $(1)\Rightarrow (2)$ As $R$ is
commutative, every full element in $R$ is invertible. Therefore
for any $a\in R$, there exists an element $e\in R$ such that
$a-e\in U(R), e^2-e\in J(R)$.

$(2)\Rightarrow (3)$ For any $a\in R$, there exist $f\in R$ and
$u\in U(R)$ such that $a=f+u$ and $f-f^2\in J(R)$. Let
$e=(1-f)+(f-f^2)u^{-1}$. Then $e-e^2\equiv f-f^2\equiv 0
\big(~mod~J(R)\big)$. Hence, $e-e^2\in J(R)$. Further,
$a-e=f+u-1+f-(f-f^2)u^{-1}=(u-2fu-u^2+f-f^2)(-u^{-1})=(a-a^2)(-u^{-1})\in
(a-a^2)R$.

$(3)\Rightarrow (4)$ For any $a\in R$, there exists $e\in R$ such
that $a-e\in (a-a^2)R, e-e^2\in J(R)$. Write $a-e=(a-a^2)s$ for
some $s\in R$. Then $e=a\big(1-(1-a)as\big)\in aR$. Further,
$1-e=(1-a)(1+as)\in (1-a)R$, as required.

$(4)\Rightarrow (1)$ For any $a\in R$, there exists $e\in R$ such
that $e\in aR, 1-e\in (1-a)R, e-e^2\in J(R)$. Write $e=as$ and
$1-e=(1-a)x$ for some $s,x\in R$. Set $f=1-e$. Then
$f-f^2=e-e^2\in J(R)$. We may assume that $se=s$ and $x(1-e)=x$.
It is easy to verify that
$(a-f)(s-x)=as-ax-fs+fx=e-ax-(1-e)es+fx=e+(f-a)x=e+(f-a)fx=e+f(1-a)x=e+f=1$.
Therefore $a=f+(a-f)$ with $a-f\in U(R)$, as asserted.\hfill$\Box$

\vskip4mm \no {\bf Corollary 5.2.}\ \ {\it Let $R$ be a
commutative ring. Then the following are equivalent:}
\vspace{-.5mm}
\begin{enumerate}
\item [(1)]{\it $R$ is feckly clean .}
\vspace{-.5mm}
\item [(2)]{\it For any $a\in R$, there exist $n\in {\Bbb N}, e\in R$ and $u\in U(R)$ such that
$a^n=e+u$ and $e-e^2\in J(R)$.}
\end{enumerate}
\vspace{-.5mm}  {\it Proof.}\ \ $(1)\Rightarrow (2)$ It is trivial
by choosing $n=1$.

$(2)\Rightarrow (1)$ For any $a\in R$, there exist $n\in {\Bbb N},
f\in R$ and $u\in U(R)$ such that $a^n=f+u$ and $f-f^2\in J(R)$.
Let $e=(1-f)+(f-f^2)u^{-1}$. Then $e-e^2\in J(R)$. One easily
checks that
$a^n-e=f+u-1+f-(f-f^2)u^{-1}=(u-2fu-u^2+f-f^2)(-u^{-1})=(a^n-a^{2n})(-u^{-1})$.
This implies that $e\in aR$. Further,
$1-e=(1-a^n)+(a^n-a^{2n})(-u^{-1})\in (1-a)R$. Therefore $R$ is
feckly clean by Theorem 5.1.\hfill$\Box$

\vskip4mm \no {\bf Corollary 5.3.}\ \ {\it Let $R$ be a
commutative ring. Then the following are equivalent:}
\vspace{-.5mm}
\begin{enumerate}
\item [(1)]{\it $R$ is feckly clean .}
\vspace{-.5mm}
\item [(2)]{\it For any $a\in R$, there exist $e\in R$ and $u\in U(R)$ such that $e=aue, 1-e=(1-a)(-u)(1-e)$ and $e-e^2\in J(R)$.}
\end{enumerate}
\vspace{-.5mm}  {\it Proof.}\ \ $(1)\Rightarrow (2)$ For any $a\in
R$, we have $a=f+v,f-f^2\in J(R), v\in U(R)$. Thus,
$a(1-f)=v(1-f)$; hence, $1-f=v^{-1}a(1-f)$. As $1-a=1-f-v$, we
have $f=(-v^{-1})(1-a)f$. Set $e=1-f$ and $u=v^{-1}$. Then
$e-e^2\in J(R)$. Moreover, $e=aue$ and $1-e=(1-a)(-u)(1-e)$, as
required.

$(2)\Rightarrow (1)$ For any $a\in R$, there exist $e\in R$ and
$u\in U(R)$ such that $e=aue, 1-e=(1-a)(-u)(1-e)$ and $e-e^2\in
J(R)$. Then $e\in aR, 1-e\in (1-a)R$. According to Theorem 5.1,
$R$ is feckly clean .\hfill$\Box$

\vskip4mm A Hausdorff space is a topological space in which
distinct points $x$ and $y$ have disjoint open sets $U$ and $V$.

\vskip4mm \no {\bf Theorem 5.4.}\ \ {\it Let $R$ be a commutative
ring. Then the following are equivalent:} \vspace{-.5mm}
\begin{enumerate}
\item [(1)]{\it $R$ is feckly clean .}
\vspace{-.5mm}
\item [(2)]{\it $Max(R)$ is a strongly zero dimensional
Hausdorff space.} \vspace{-.5mm}
\item [(3)]{\it $Max(R)$ is strongly zero dimensional and $a+b=1$ in $R$ implies that $(1+ar)(1+bs)\in J(R)$
for some $r,s\in R.$} \vspace{-.5mm}
\item [(4)]{\it $a+b =1$ in $R$ implies that $1+ar\in eR,
1+bs\in (1-e)R$ where $e-e^2\in J(R)$.}
\end{enumerate}
\vspace{-.5mm}  {\it Proof.}\ \ $(1)\Leftrightarrow (2)$ In view
of [7, Theorem 1.3], $Max(R)$ is an Hausdorff space if and only if
every prime ideal of $R$ containing $J(R)$ is contained in a
unique maximal ideal, and proving $(2)$ by Theorem 3.5.

$(1)\Rightarrow (3)$ In view of Theorem 3.5, $Max(R)$ is strongly
zero dimensional. Given $a+b =1$ in $R$, there exists $u\in U(R)$
such that $a=e+u$ where $e-e^2\in J(R)$. It is easy to verify that
$$\begin{array}{lll}
&\big(1+\frac{a}{e-a}\big)\big(1-\frac{b}{e-a}\big)\\
=&1-\frac{b}{e-a}+\frac{a}{e-a}-\frac{ab}{(e-a)^2}\\
=&\frac{e^2-e}{(e-a)^2}\\
\in &J(R).
\end{array}$$ Choose $r=\frac{1}{e-a}$ and $t=-\frac{1}{e-a}$.
Then $(1+ar)(1+bs)\in J(R)$.

$(3)\Rightarrow (1)$ Let $P$ be a prime ideal of $R$ such that
$J(R)\subseteq P$ and $P$ is contained in two distinct maximal
ideals $M$ and $N$. Write $a+b=1$ with $a\in M, b\in N$. By
hypothesis, there are some $r,s\in R$ such that $(1+ar)(1+bs)\in
J(R)$. As $J(R)\subseteq P$, we see that either $1+ar\in P$ or
$1+bs\in P$. If $1+ar\in P$, then $1=(1+ar)-ar\in M$. If $1+bs\in
P$, then $1=(1+bs)-bs\in N$. In any case, we get a contradiction.
Therefore every prime ideal of $R$ containing $J(R)$ is contained
in exactly one maximal ideal. According to Theorem 3.5, $R$ is
feckly clean.

$(1)\Rightarrow (4)$ Suppose that $R$ is feckly clean and $a+b=1$
in $R$. In view of Theorem 5.1, there exists an element $f\in R$
such that $f\in aR, 1-f\in (1-a)R$ and $f^2-f\in J(R)$. Write
$1-f=-(1-a)s$ and $f=-ar$ for some $r,s\in R$. Then $1+bs=f$ and
$1+ar=1-f$. Set $e=1-f$. Then $e^2-e=f^2-f\in J(R)$. Further,
$1+ar\in eR$ and $1+bs\in (1-e)R$, as required.

$(4)\Rightarrow (1)$ is trivial by Corollary 2.7.\hfill$\Box$

\vskip4mm \no {\bf Corollary 5.5.}\ \ {\it A commutative ring $R$
is feckly clean if and only if $Max(R)$ is strongly
zero-dimensional and $J$-$spec(R)$ is normal.} \vskip2mm\no {\it
Proof.}\ \ In view of [7, Theorem 1.3], $J$-$spec(R)$ is normal if
and only if $Max(R)$ is a Hausdorff space, and so the result
follows by Theorem 5.4.\hfill$\Box$

\vskip4mm Recall that a commutative ring $R$ is a $pm$ (Gelfand)
ring provided that each prime ideal of $R$ is contained in exactly
one maximal ideal. In [1, Corollary 4], Anderson and Camillo
proved that every commutative clean ring is always a $pm$ ring.

\vskip4mm \no {\bf Theorem 5.6.}\ \ {\it Let $R$ be a commutative
ring. Then the following are equivalent:} \vspace{-.5mm}
\begin{enumerate}
\item [(1)]{\it $R$ is clean.}
\vspace{-.5mm}
\item [(2)]{\it $R$ is feckly clean and $R$ is a $pm$ ring.}
\vspace{-.5mm}
\item [(3)]{\it $R$ is feckly clean ring in which $a+b
=1$ implies that $(1+ar)(1+bs)=0$ for some $r,s\in R.$}
\end{enumerate} \vspace{-.5mm} {\it Proof.}\ \ $(1)\Rightarrow
(2)$ Clearly, $R$ is feckly clean. According to [1, Corollary 4],
$R$ is a $pm$ ring.

$(2)\Rightarrow (1)$ For any disjoint closed sets $A$ and $B$ of
$Max(R)$, it follows from Theorem 2.3 that there exists an element
$e\in R$ such that $A\subseteq V(e), B\subseteq V(1-e)$ and
$eR(1-e)\subseteq J(R)$. As $eR(1-e)\subseteq J(R)$, we get
$Max(R)=V(e)\bigsqcup V(1-e)$; hence that $V(1-e)$ and $V(e)$ are
clopen sets of $Max(R)$. Therefore $Max(R)$ is strongly
zero-dimensional. According to [4, Corollary 17.1.14], $R$ is
strongly clean.

$(2)\Leftrightarrow (3)$ According to [5, Theorem 4.4], $R$ is a
pm ring if and only if $a+b =1$ implies that $(1+ar)(1+bs)=0$ for
some $r,s\in R$, and so the result follows.\hfill$\Box$

\vskip4mm \no  {\bf Corollary 5.7.}\ \ {\it Let $X$ be a
topological space. Then the following are equivalent:}
\vspace{-.5mm}
\begin{enumerate}
\item [(1)]{\it $C(X)$ is clean.}
\vspace{-.5mm}
\item [(2)]{\it $C(X)$ is feckly clean.}
\end{enumerate} \vspace{-.5mm} {\it Proof.}\ \ In view of [3],
$C(X)$ is a pm ring. Thus, the proof is true by Theorem
5.6.\hfill$\Box$

\vskip4mm \no  {\bf Corollary 5.8.}\ \ {\it Let $R$ be a
commutative ring. Then the following are equivalent:}
\vspace{-.5mm}
\begin{enumerate}
\item [(1)]{\it $R$ is $\pi$-regular.}
\vspace{-.5mm}
\item [(2)]{\it $R$ is feckly clean ring in which every prime ideal is maximal.}
\end{enumerate} \vspace{-.5mm} {\it Proof.}\ \ $(1)\Rightarrow
(2)$ In view of Theorem 4.1, $R$ is feckly clean. According to [6,
Corollary 2.8], every prime ideal of $R$ is maximal.

$(2)\Rightarrow (1)$ In light of Theorem 3.5, $Max(R)$ is strongly
zero-dimensional. Therefore we complete the proof by [6, Theorem
2.3].\hfill$\Box$

\vskip10mm \bc{\small ACKNOWLEDGEMENTS}\ec \vskip4mm The authors
are grateful to the referee for his/her helpful suggestions which
correct many errors in the first version and lead the new one more
clear. This research was supported by the Natural Science
Foundation of Zhejiang Province (LY13A010019) and the Scientific
and Technological Research Council of Turkey (2221 Visiting
Scientists Fellowship Programme).

\end{document}